\documentclass[12 pt]{article}

\makeatletter
\def\blfootnote{\gdef\@thefnmark{}\@footnotetext}
\makeatother

\usepackage{amssymb} 
\usepackage{amsmath} 
\usepackage{latexsym}
\usepackage{graphicx}

\usepackage{enumerate}

\usepackage{tikz}
\usetikzlibrary{calc}
\usepackage{pgf,tikz}
\usetikzlibrary{arrows}
\usepackage{hyperref}
\usepackage{amscd}
\usepackage{relsize}
\usepackage{verbatim}

\usepackage{tikz-cd}

\usepackage{ color, amsmath, amssymb, amsfonts, amscd, graphicx,latexsym, hyperref}
\usepackage{caption, subcaption}
\usepackage[pdf]{pstricks}

\usepackage[affil-it]{authblk}

\usepackage{amsthm}

\newcommand{\R}{\mathbf R}

\newtheorem{theorem}{Theorem}

\newtheorem{proposition}{Proposition}

\newtheorem{remark}{Remark}

\newtheorem{lemma}{Lemma}

\newtheorem{definition}{Definition}

\providecommand{\keywords}[1]
{
	\small	
	\textbf{\textit{Keywords---}} #1
}

\title{Classification of punctures on complete flat surfaces}
\author{ \.{I}smail Sa\u{g}lam \thanks{Electronic address: \texttt{isaglamtrfr@gmail.com}  }}
\affil{Adana Alparslan T\"{u}rke\c{s} Science and Technology University  }

\date{}

\begin{document}

\maketitle
\blfootnote{\textup{2000} \textit{Mathematics Subject Classification}:
51F99, 57M50}
\begin{abstract}
We investigate the behavior of a complete flat metric on a surface near a puncture. We call a puncture on  a flat surface regular if it has a neighborhood which is isometric to that of a point at infinity of a cone. We prove that there are punctures which are not regular if and only if the curvature at the puncture is $4\pi$. 
\end{abstract}


\thanks{}
\section{Introduction}
\thanks{}
Flat surfaces are obtained by gluing Euclidean triangles along their edges appropriately. They appear in several areas of mathematics and physics. For example, they are studied in dynamics of billiard tables. It is known that each rational polygon can be covered by a flat surface with trivial holonomy group. Such surfaces are called translation surfaces and have been studied  extensively \cite{dynamics}. Also, these surfaces are quite useful in Teichm\"{u}ller theory. Together with quadratic differentials, they are used in the proofs of Teichm\"{u}ller's theorems \cite{teichmuller}. They appear in quantum gravity and topological quantum field theory as well.  See \cite{gravity} and \cite{quantum-field}.

These surfaces are  interesting for their own sake. For example, Thurston obtained complex hyperbolic orbifolds from moduli spaces of certain flat spheres \cite{thurston}. Following Thurston, Bavard and Ghys obtained real hyperbolic orbifolds from the moduli spaces of certain polygons in the plane \cite{bavard-ghys}. Troyonov introduced certain geometric structures on Teichm\"{u}ller spaces by considering the moduli spaces of flat surfaces with prescribed curvature data \cite{Tro-handbook}.

Compact flat surfaces are examples of length spaces. There is a length minimizing geodesic between any two points of such a surface. 
\cite{Gromov}, \cite{dimitri}. They can be triangulated with finitely many triangles.  In addition, Gauss-Bonnet formula holds for these surfaces.
See \cite{Tro-enseign},  \cite{Tro-compact}.

Flat surfaces with \textit{regular} punctures have been  studied in \cite{Tro-open}. By a regular puncture on a flat surface, we mean a puncture which has a neighborhood isometric to that of the point at infinity of a cone. Flat surfaces with possibly irregular punctures have been studied in \cite{saglam}. We now state the main results of \cite{saglam}. Let $\bar{S}$ be a complete flat surface.

\begin{enumerate}
\item 
$\bar{S}$ can be triangulated with finitely many types of triangles. 
\item
 Gauss-Bonnet formula holds for $\bar{S}$. 
\item
Each loop on $\bar{S}$ has a geodesic representative in its free homotopy class.
	\end{enumerate}

Our objective is to understand the behavior of a complete flat metric near a puncture. More precisely,  we classify complete flat metrics on a disk with a puncture up to the  modification equivalence, where two flat metrics on a disk are equivalent if they are "same" on a neighborhood of the puncture. Now we state the main results of the present paper.

Let $\bar{S}$ be a flat surface. 

\begin{itemize}
	\item 
	If the curvature at a punctured interior point  equals to $4\pi$, then there are uncountably many 
	modification-equivalence classes of complete flat metrics near the puncture. 
	\item
	If the curvature at a punctured interior point  is not equal to $4\pi$, then any complete flat metric is modification-equivalent to a cone near the puncture.
	
\end{itemize}

\subsection{Doubly labeled surfaces}
\label{notation}
In this paper, we use the notation in \cite{saglam}. A doubly labeled surface is a compact surface together with labeled points.

\begin{definition}
	Let $S$ be a connected compact topological  surface perhaps with boundary $B$.
	Let $\mathfrak{l}, \mathfrak{p},\mathfrak{l'},\mathfrak{p'}$ be finite disjoint subsets of $S$ so that 
	\begin{itemize}
		\item
		$\mathfrak{l}$ and $\mathfrak{p}$ are  subsets of the \textit{interior} of $S$,
		\item
		$\mathfrak{p'}, \mathfrak{l'}$ are  subsets of $B$.
		
	\end{itemize}

	An element in $\mathfrak{l}$ will be called  a \textit{labeled interior}  point.  An element in $\mathfrak{p}
	$ will be called a punctured interior point. Other points in the interior of $S$ are  called \textit{ordinary interior points}. An element in $B$ will be called  a boundary point. An element
	in $\mathfrak{l'}$ will be called a \textit{labeled boundary} point. An element in $\mathfrak{p'}$
	will be called  \textit{punctured boundary} point.
	Other points in the boundary will be called \textit{ordinary boundary points}. 
	A \textit{doubly labeled surface}, shortly DL surface, is the tuple 
	$$ (S,B, \mathfrak{l},\mathfrak{p}, \mathfrak{l'}, \mathfrak{p'} )$$
\end{definition}

Also we will use the following notation:
\begin{enumerate}
	\item
	$S_B=S-B$.
	\item
	$S_{\frak{l}}=S-\frak{l}$
	\item
	$S_{B,\frak{l}}=S-(B\cup \frak{l})$
	\item
	\dots
\end{enumerate}

We will denote a doubly labeled surface $ (S,B, \mathfrak{l},\mathfrak{p},\mathfrak{l'}, \mathfrak{p'} ) $ as $S^L$. Underlying compact surface of $S^L$ will simply be denoted by $S$. 
Note that DL surfaces can be considered as punctured surfaces with puncture set  $\frak{p} \cup \mathfrak{p'}$. Indeed, $S_{\frak{p,p'}}$ is the punctured surface that we consider. 
We point out that the punctured and labeled points may lie in the boundary.

A cone having angle $\theta>0$, or equivalently  curvature $\kappa =2\pi-\theta$, is the set
\begin{align}
\{(r, \psi): r \in {\R}^{\geq 0}, \psi \in {\R}/\theta \mathbb{Z}\}
\end{align}
with the metric
\begin{align}
\mu=dr^2+r^2d\psi^2.
\end{align}

 We can consider a cone  as a  DL   sphere with one punctured and one labeled interior point. The point $(0,0)$ is called vertex or the origin of the cone. 
We will denote it  by $v_{0}$. Let $\theta(v_0)=\theta$ and $\kappa(v_0)=2\pi - \theta$. Note that we may  talk about the \textit{point at infinity}  or the punctured point.
We shall denote this point by $v_{\infty}$.

\begin{definition}
	Consider a cone with angle $\theta >0$.
	\begin{enumerate}
		\item
		$\kappa{(v_\infty)}=2\pi+\theta$ is called the curvature at $v_{\infty}$.
		\item
		$\theta{(v_\infty)}=-\theta$ is called the angle at $v_{\infty}$.
		
	\end{enumerate}
\end{definition}

\noindent  We will denote 
a cone with angle $\theta$ by $C_{\theta}$.

\begin{definition}
	A cylinder of width $r$, $C_{0r}$, is a metric space obtained by identifying edges of  an infinite strip in the Euclidean plane  having  width $r$ through \textit{opposite} points.
\end{definition}

\noindent Observe that a cylinder can be considered as  DL sphere with two punctured points. By convention, the angles at these punctures are $0$. We  also call a cylinder as a cone of angle $0$. Also, again by convention, the curvature at each of the punctured points,  is $2\pi$.

\begin{definition}
	A (flat) cone metric on a DL surface $S^L$ is a metric on $S_{\frak{p,p'}}$ so that each point $x$ in  $S_{\frak{p,p'}}$ has a neighborhood isometric to a neighborhood of the apex of the cone $C_{\theta}=C_{\theta_x}$ or a section of a cone $V_{\theta}=V_{\theta_x}$, and
	\begin{enumerate}
		\item
		$\frak{l}=\{y \in S_{\frak{p},B}: \ \theta_y \neq 2\pi \}$,
		\item
		$\frak{l'}= \{ y \in B-\frak{p'}: \ \theta_y \neq \pi   \}$.
	\end{enumerate}
	Angle at $x$, $\theta(x)$ is defined to be $\theta_x$. If $x \in S_{\frak{p},B}$, then 
	the	curvature at $x$, $\kappa(x)$, is defined as $2\pi - \theta(x) $. If $x \in B-\frak{p'}$, then the curvature is  $\kappa(x)=\pi -\theta(x)$. $x$ is called 
	singular if $\kappa(x)\neq 0$. Otherwise it is called non-singular.
	
\end{definition}
\noindent Note that the conditions {\it 1.} and {\it 2.} assure that the set of singular points and $\frak{l} \cup \frak{l'}$ are same.
A flat surface with a flat cone metric is called a flat DL surface.

\begin{definition}
	A punctured interior point  on a flat DL surface is called regular if it has a neighborhood  isometric to a neighborhood of the point at infinity of a cone. Otherwise, it is called irregular. 
\end{definition}

An example of a flat DL surface with an irregular punctured interior point is given in \cite{saglam}. Also, the curvature at a punctured interior (or boundary) point of a flat DL surface was defined in \cite{saglam}. For a DL surface, the following theorem holds. See \cite{saglam} for a proof.

 \begin{theorem}[Gauss-Bonnet formula]
 	Let $S^L$ be a complete flat DL surface. The following formula holds: 
 	
 	\begin{align}
 	\sum_{x\in S}\kappa(x)= 2\pi\chi(S),
 	\end{align}
 	where $\chi(S)$ is the Euler characteristics of $S$.
\end{theorem}

\subsection{Modification }
\label{modification1}

If we have a complete flat metric on a DL surface with boundary, then we can cut a triangle having one edge incident to the boundary to get another complete flat metric. Note that the behavior of the metric near the punctures remains unchanged after this operation. By a modification of flat DL surface, we mean the surface obtained by removing finitely many triangles which are incident to the boundary. 

We denote a closed  disk with one punctured interior point by $D^L$.  We assume that $D^L$ has no labeled interior points and punctured boundary points. As usual, we denote the underlying closed disk by $D$. Let $\bar{D}^L$ be a closed disk with one punctured boundary point. We assume that $\bar{D}^L$ has no labeled interior points or punctured interior points. We denote the underlying closed disk by $\bar{D}$.

\begin{proposition}[\cite{saglam}]
	\label{mod2}
	Each complete cone metric on $D^L$ can be modified so that resulting disk
	does not have any points with positive curvature on its boundary.

\end{proposition}

\section{ Flat  metrics on a disk with one punctured point}
\label{modification}

\subsection{Modification Equivalence}
Recall that a modification 
of a flat metric on $D^L$ is a flat metric obtained by successively cutting Euclidean triangles which are incident to its boundary. In this section, we will classify  complete flat metrics on $D^L$ up to modification equivalence. Two flat complete metrics $\mu$ and $\eta$ are called the modification equivalent if they can be modified so that there is an orientation preserving isometry between the resulting complete flat disks.  Lemma \ref{mod2}  implies that any flat complete   metric on $D^L$ contains a metric with non-positive curvature data in its equivalence class. Thus, from now on, by a complete flat metric on $D^L$, we mean a metric with non-positive curvature data. We will also study the case of the disc with one punctured boundary point. Note that one can define modification equivalence for flat complete metrics on $\bar{D}^L$ in a similar way.

\subsection{Reducing number of singular points}
Let $\mu$ be a complete  cone metric on $D^L$ and $\frak{K}$ be the total curvature on its boundary. If $\mathfrak{K}=0$ then $D^L$ is modification isometric to a half-cylinder. Therefore we need to consider the case $\mathfrak{K}<0$. 
\noindent We start with a simple fact.

\begin{lemma}
	Let $a_1,\dots,a_k$ be real numbers so that not all of them are equal. Let 
	\begin{equation}
	a=\frac{\sum_{i=1}^k a_i}{k}
	\end{equation}
	be their avarage. There exists i so that 
	
	\begin{itemize}
		\item 
		$a_i\geq a$ and $a_{i+1} < a$, or
		\item 
		$a_i < a$ and $a_{i+1} \geq a$,
	\end{itemize}
	where $a_{k+1}=a_1$.
	

	
\end{lemma}

\begin{figure}
	\begin{center}
		\includegraphics[scale=0.5]{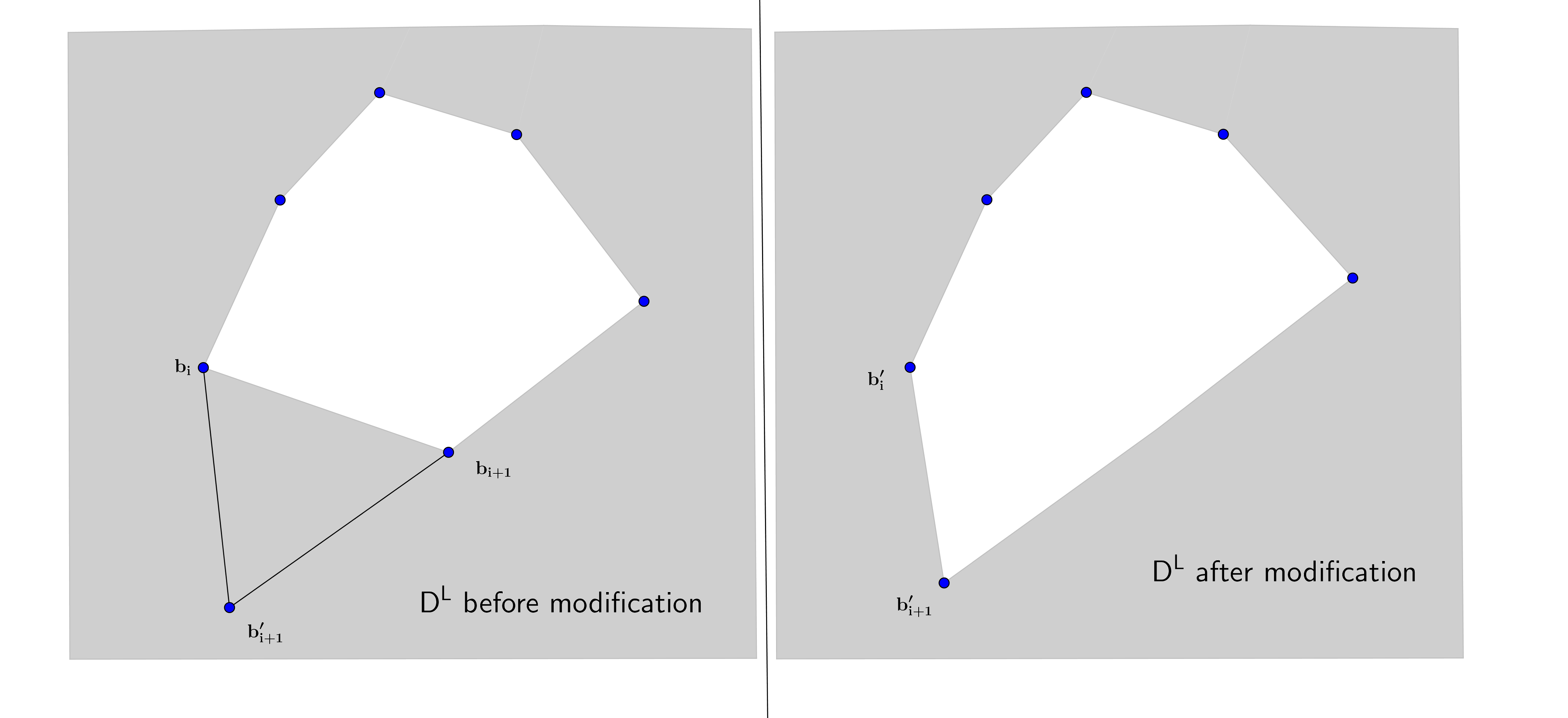}
	\end{center}
	\caption{After the modification, curvature at 
		the $i$-th vertex $b_i$ becomes $\frac{\frak{K}}{n}$.}
	\label{fig.cut1}
\end{figure}

\begin{lemma}
	\label{homojen}
	Assume that $D^L$ has $k$ labeled boundary points. $D^L$ can be modified 
	so that resulting metric has also $k$ singular points and these points have same curvature.
	\begin{proof}
		First label singular points as $b_1, \dots, b_k$ so that they are in a cyclic order on the boundary. Assume that not all of the curvatures are equal. By above lemma, there exists $i$ so that either
		$$\kappa_i<\frac{\frak{K}}{k}\ \textup{and} \ \kappa_{i+1} \geq\frac{\frak{K}}{k}.$$
		or,
		$$\kappa_i \geq \frac{\frak{K}}{k}\ \textup{and} \ \kappa_{i+1} < \frac{\frak{K}}{k},$$

		\noindent Let us consider the first case. Remove the triangle with vertices 
		$b_i,b_{i+1},b_{i+1}'$ having angles $\lvert \kappa_i \rvert - \frac{\lvert \frak{K} \rvert}{k}, \ \lvert \kappa_{i+1} \rvert, \ \pi + \frac{\lvert \frak{K} \rvert}{k}-\lvert \kappa_i \rvert - \lvert \kappa_{i+1}\rvert$ from $D^L$, respectively. See Figure \ref{fig.cut1}.
		\begin{itemize}
			
			\item
			
			if $\kappa_{i+1}>\frac{\frak{K}}{k}$, then resulting 
			metric has less singular points having curvature not equal to
			$\frac{\frak{K}}{n}$, 
			\item
			if $\kappa_{i+1}=\frac{\frak{K}}{k}$, then $i$-th vertex has curvature $\kappa_{i+1}=\frac{\frak{K}}{k}$ and $i+1$-th
			vertex has curvature $\kappa_i$ in the resulting disc.
		\end{itemize}


		Observe that if $\kappa_i \geq \frac{\frak{K}}{k}\ \textup{and} \ \kappa_{i+1} < \frac{\frak{K}}{k}$, there is a similar cutting operation with the following properties:
		
		\begin{itemize}
			\item
			if $\kappa_{i}>\frac{\frak{K}}{k}$ then resulting metric resulting metric has less singular points having curvature not equal to 
			$\frac{\frak{K}}{k}$.
			\item
			if $\kappa_{i}=\frac{\frak{K}}{k}$, then  $i+1$-th vertex has curvature $\frac{\frak{K}}{k}$ and $i$-th
			vertex has curvature $\kappa_{i+1}$ in the resulting 
			disc.
		\end{itemize}

		Apply the following algorithm to $D^L$ repeatedly.
		\begin{enumerate}
			\item
			If there is an index $i$ so that either 
			$$\kappa_i > \frac{\frak{K}}{k}\ \textup{and} \ \kappa_{i+1} < \frac{\frak{K}}{k},$$
			or,
			$$\kappa_i<\frac{\frak{K}}{k}\ \textup{and} \ \kappa_{i+1} >\frac{\frak{K}}{k}.$$
			apply cutting operation described above. Observe that the number of
			singular points having curvature which is not equal to $\frac{\frak{K}}{k}$
			decreases after this operation. Do this repeadetly so that 
			there are no index $i$ satifying any of the properties above.
			\item
			After the first step, if all curvatures of the singular points of the resulting metric are equal, then we are done. If this is not the case, we can permute curvatures to  get a modification for which 
			there exists $i$ so that 
			$$\kappa_i > \frac{\frak{K}}{k}\ \textup{and} \ \kappa_{i+1} < \frac{\frak{K}}{k},$$
			or,
			$$\kappa_i<\frac{\frak{K}}{k}\ \textup{and} \ \kappa_{i+1} >\frac{\frak{K}}{k}.$$
			\item
			Apply the first step to the resulting metric.
		\end{enumerate}
		Since the number of singular points having curvature not equal 
		to $\frac{\frak{K}}{k}$ decreases at each \textit{run} of the algorithm, we get a cone metric of desired type in finitely many steps.
		
	\end{proof}
	
\end{lemma}

\begin{figure}
	\begin{center}
		\includegraphics[scale=0.4]{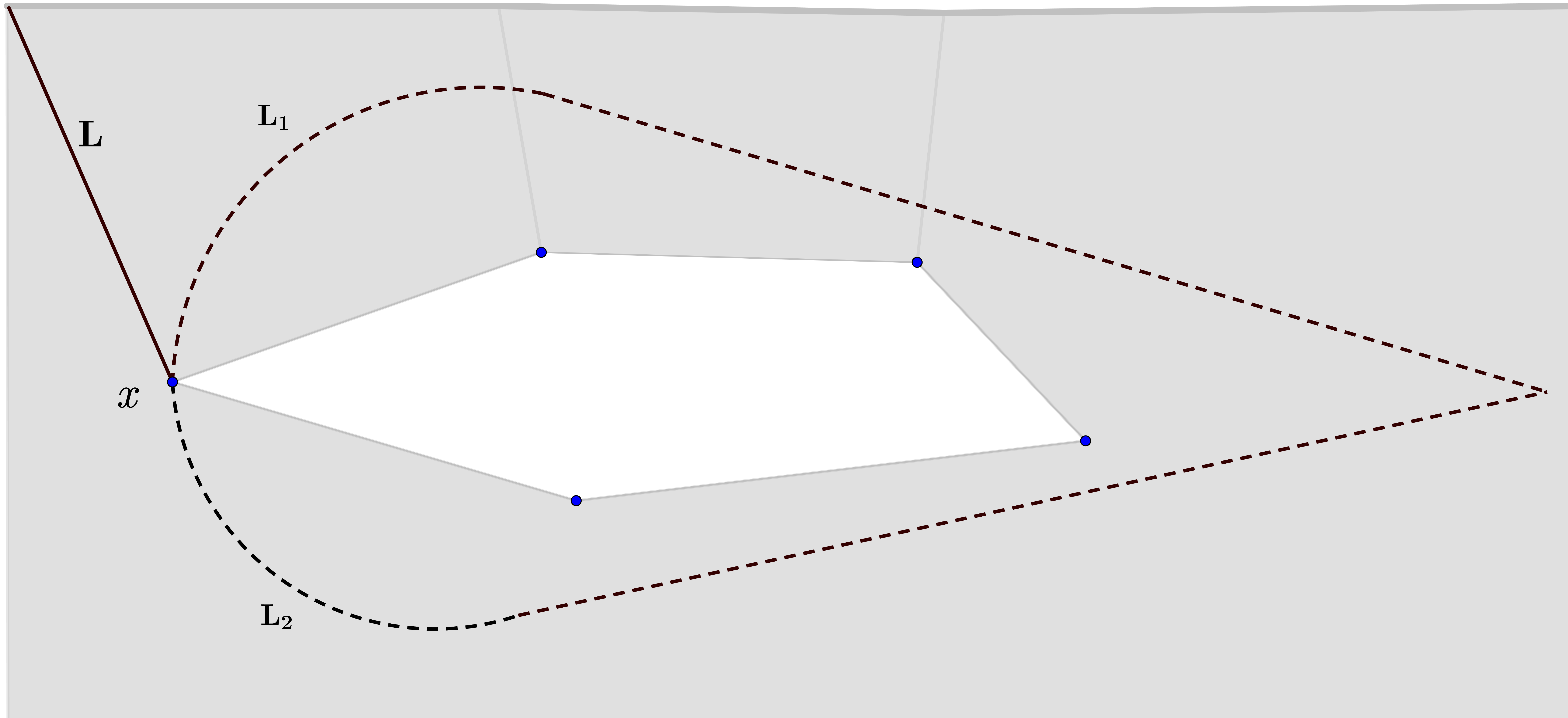}
	\end{center}
	\caption{Cut the annulus bounded by $L_1,L_2$  and the boundary of $D^L$ to obtain a cone metric with one singular point having curvature $\frak{K}$.}
	\label{reduce-to-vertex}
\end{figure}

\begin{lemma}
	\label{base-reduce}
	Assume that $0\leq\lvert \frak{K} \rvert < \pi$. Let $x$ be a boundary point on  $D^L$ and $L$ be a half-line originating from $x$ and directing toward its interior.  This disc can be modified so that resulting metric has at most  \textit{one} singular point and no points except $x$ on $L$ is in the  triangles removed during the modification.
	\begin{proof}
		The case $\frak{K}=0$ is trivial since there are no singular points for this case.
		Desired modification for the cone metrics having singular points is described in Figure \ref{reduce-to-vertex}.
		Draw two half lines $L_1,L_2$ originating from $x$ making an angle of $\pi$ with each other. Gauss-Bonnet formula implies that $L_1$ and $L_2$ intersect only at 2 points, one of which is $x$. These 
		lines, together with the boundary of $S^L$, bound a compact region with polygonal geodesic boundary. We can obtain desired modification by removing this region. 
	\end{proof}
	
\end{lemma}

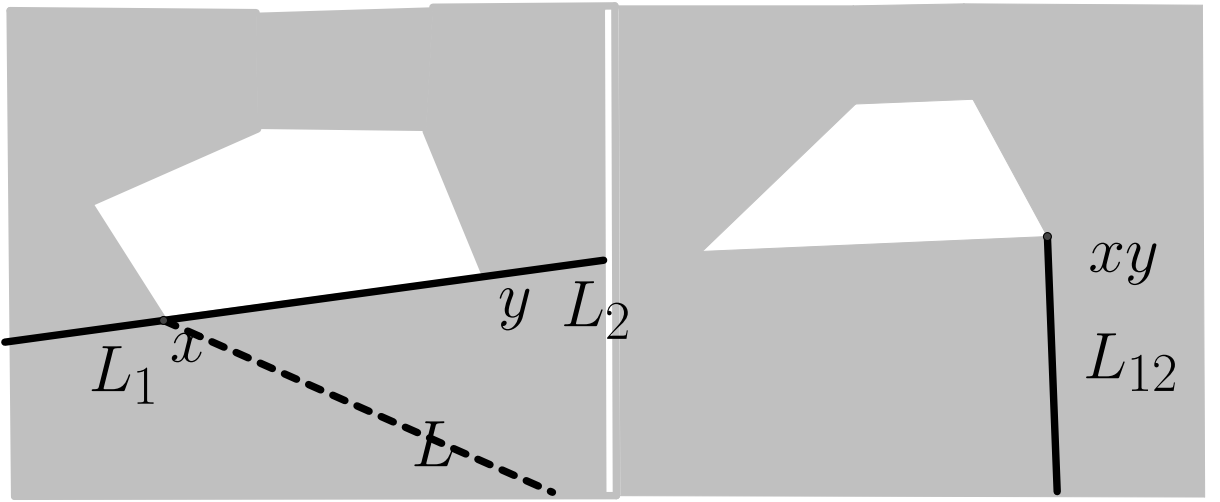
\begin{figure}
	\begin{center}
		\definecolor{uuuuuu}{rgb}{0.26666666666666666,0.26666666666666666,0.26666666666666666}
		\definecolor{cqcqcq}{rgb}{0.7529411764705882,0.7529411764705882,0.7529411764705882}
		\definecolor{ffffff}{rgb}{1.,1.,1.}
		\begin{tikzpicture}[line cap=round,line join=round,>=triangle 45,x=0.5cm,y=0.45cm]
		
		\fill[line width=2.8pt,color=cqcqcq,fill=cqcqcq,fill opacity=0.5] (2.42,2.66) -- (2.38,6.1) -- (-4.16,6.16) -- (-4.04,-8.2) -- (11.96,-8.18) -- (11.92,6.3) -- (7.08,6.26) -- (6.9,2.6) -- (8.5,-1.72) -- (-0.08,-3.) -- (-2.06,0.46) -- cycle;
		\fill[color=cqcqcq,fill=cqcqcq,fill opacity=0.5] (2.38,6.1) -- (7.08,6.26) -- (6.9,2.6) -- (2.42,2.66) -- cycle;
		\fill[color=cqcqcq,fill=cqcqcq,fill opacity=0.5] (18.32,3.4) -- (18.26,6.3) -- (11.92,6.3) -- (11.96,-8.18) -- (27.6,-8.22) -- (27.54,6.32) -- (21.2,6.36) -- (21.44,3.54) -- (23.42,-0.52) -- (14.24,-0.96) -- cycle;
		\fill[color=cqcqcq,fill=cqcqcq,fill opacity=0.5] (18.26,6.3) -- (18.32,3.4) -- (21.44,3.54) -- (21.2,6.36) -- cycle;
		\draw [line width=5.2pt,color=ffffff] (11.84,6.14)-- (11.88,-8.38);
		\draw [line width=2.8pt,color=cqcqcq] (2.38,6.1)-- (-4.16,6.16);
		\draw [line width=2.8pt,color=cqcqcq] (-4.16,6.16)-- (-4.04,-8.2);
		\draw [line width=2.8pt,color=cqcqcq] (-4.04,-8.2)-- (11.96,-8.18);
		\draw [line width=2.8pt,color=cqcqcq] (11.96,-8.18)-- (11.92,6.3);
		\draw [line width=2.8pt,color=cqcqcq] (11.92,6.3)-- (7.08,6.26);
		\draw [line width=2.8pt,color=cqcqcq] (6.9,2.6)-- (8.5,-1.72);
		\draw [line width=2.8pt,color=cqcqcq] (8.5,-1.72)-- (-0.08,-3.);
		\draw [line width=2.8pt,color=cqcqcq] (-0.08,-3.)-- (-2.06,0.46);
		\draw [line width=2.8pt,color=cqcqcq] (-2.06,0.46)-- (2.42,2.66);
		\draw [line width=2.8pt,domain=-4.3:11.62] plot(\x,{(--34.9576-1.78*\x)/-11.7});
		\draw [color=cqcqcq] (7.08,6.26)-- (6.9,2.6);
		\draw [color=cqcqcq] (2.42,2.66)-- (2.38,6.1);
		\draw [color=cqcqcq] (18.32,3.4)-- (18.26,6.3);
		\draw [color=cqcqcq] (18.26,6.3)-- (11.92,6.3);
		\draw [color=cqcqcq] (11.92,6.3)-- (11.96,-8.18);
		\draw [color=cqcqcq] (11.96,-8.18)-- (27.6,-8.22);
		\draw [color=cqcqcq] (27.6,-8.22)-- (27.54,6.32);
		\draw [color=cqcqcq] (27.54,6.32)-- (21.2,6.36);
		\draw [color=cqcqcq] (21.2,6.36)-- (21.44,3.54);
		\draw [color=cqcqcq] (21.44,3.54)-- (23.42,-0.52);
		\draw [color=cqcqcq] (23.42,-0.52)-- (14.24,-0.96);
		\draw [color=cqcqcq] (14.24,-0.96)-- (18.32,3.4);
		\draw [color=cqcqcq] (18.26,6.3)-- (18.32,3.4);
		\draw [color=cqcqcq] (18.32,3.4)-- (21.44,3.54);
		\draw [color=cqcqcq] (21.44,3.54)-- (21.2,6.36);
		\draw [color=cqcqcq] (21.2,6.36)-- (18.26,6.3);
		\draw (-0.2,-3.06) node[anchor=north west] {{\huge $x$}};
		\draw (8.52,-1.74) node[anchor=north west] {{\huge $y$}};
		\draw (-2.36,-3.44) node[anchor=north west] {{\huge $L_1$}};
		\draw (10.2,-1.52) node[anchor=north west] {{\huge $L_{2}$}};
		\draw (24.08,-3.06) node[anchor=north west] {{\huge $L_{12}$}};
		\draw (6.22,-5.66) node[anchor=north west] {{\huge $L$}};
		\draw (24.22,-0.4) node[anchor=north west] {{\huge $xy$}};
		\draw [line width=2.8pt,dash pattern=on 5pt off 5pt] (-0.08,-3.)-- (10.26,-8.08);
		\draw [line width=2.8pt] (23.42,-0.52)-- (23.68,-8.06);
		\begin{scriptsize}
		\draw [fill=uuuuuu] (-0.08,-3.) circle (1.5pt);
		\draw [fill=uuuuuu] (23.42,-0.52) circle (1.5pt);
		\end{scriptsize}
		\end{tikzpicture}
	\end{center}
	\caption{Remove  half plane determined by $L_1, L_2 $ and boundary segment $[x,y]$ and glue 
		$L_1$ and $L_2$ to get a cone metric of total curvature curvature   $\frak{K'}$, 	
		$ (n-1)\pi \leq \lvert \frak{K'} \rvert < n\pi$. Modify it, then cut it thorough $L_{12}$
		and add the half plane removed to get a desired modification of the cone metric we started.}
	\label{L1-L2}
	
\end{figure}

\begin{proposition}
	\label{L-cut}
	Assume that $n>0$ and $\frak{K} \leq 0$ so that
	
	$$(n-1)\pi \leq \lvert \frak{K} \rvert < n \pi.$$
	Also assume that the metric has $n+m, \ m\geq 0$ singular points. Let $x$ be one of them  and $L$
	be a half-line originating from $x$ and pointing interior of the disc. 
	The metric can be modified so that the resulting metric has $n$ singular points of negative curvature and the removed triangles contain no points on $L$ other than $x$.
	
	\begin{proof}
		We will prove the statement by induction on the number $n$. The base case $n=1$ is done by 
		the Lemma \ref{base-reduce}. Assume that $n\geq 1$ and 
		
		$$ n\pi \leq \lvert \frak{K} \rvert < (n+1)\pi.$$
		
		Take a cone metric on the disc having $n+1+m$ singular points. Take a boundary segment with vertices $x,y$ and half-lines $L_1, L_2$ making angles $\pi$ with the segment.  See Figure \ref{L1-L2}. Cut the half plane and glue the half-lines $L_1$ and $L_2$. By this way we get a new cone metric on the punctured disc so that the total curvature  at its boundary is $\frak{K'}=\frak{K}+\pi$, thus satisfies the inequality below:
		
		$$(n-1)\pi \leq \lvert \frak{K'} \rvert < n\pi.$$ 
		Let $L_{12}$ the half line formed by gluing $L_1$ with $L_2$.  We show the vertex obtained by $x$ and $y$ by $xy$. Observe that by induction hypothesis, we can modify this new metric so that the resulting metric has $n$ singular points and removed triangles 
		do not contain any point of $L_{12}$ except $x$. Now cut this new punctured disc together with the induced metric through $L_{12}$, and glue the half-plane 
		that we removed as in the figure. Resulting cone metric has $n+1$ singular points and is a modification of the metric we started with. Also observe that compact part removed during the modification does not intersect with half-line $L$ except at the vertex $x$.
	\end{proof}

\end{proposition}

\begin{theorem}
	\label{best-reduction}
	Every complete cone metric on $D^L$ with total boundary curvature 
	$$(n-1)\pi \leq \lvert \frak{K}\rvert < n\pi$$
	
\noindent 	can be modified so that resulting metric has $n$ singular points, and curvature of each of these singular points is $\frac{\frak{K}}{n}$. 
	
	\begin{proof}
		First modify the metric as in  Proposition \ref{L-cut} to get a metric with $n$ singular points of negative curvature. Modify this new metric as in  Lemma \ref{homojen},  to get a metric of desired type.
	\end{proof}
	
\end{theorem}

\subsection{The case $0 \leq \lvert \frak{K} \rvert < 2\pi$}

\begin{figure}
	\begin{center}
		\includegraphics[scale=0.5]{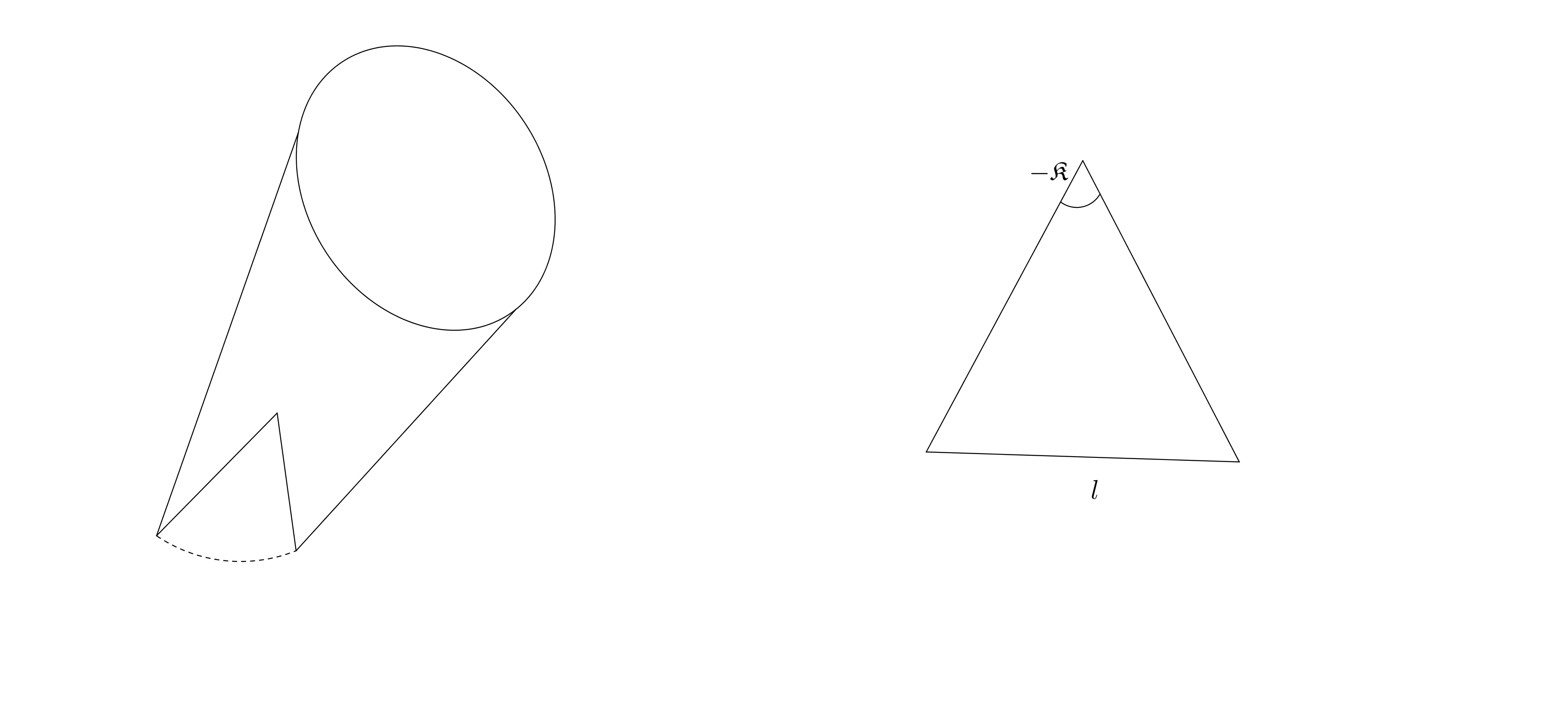}
	\end{center}
	\caption{Glue unlabeled edges of the isosceles triangle to get a flat disk. Then glue this disk with $D^L$ to get a cone.}
	\label{cone1}
\end{figure}

\begin{figure}
	\begin{center}
		\includegraphics[scale=0.50]{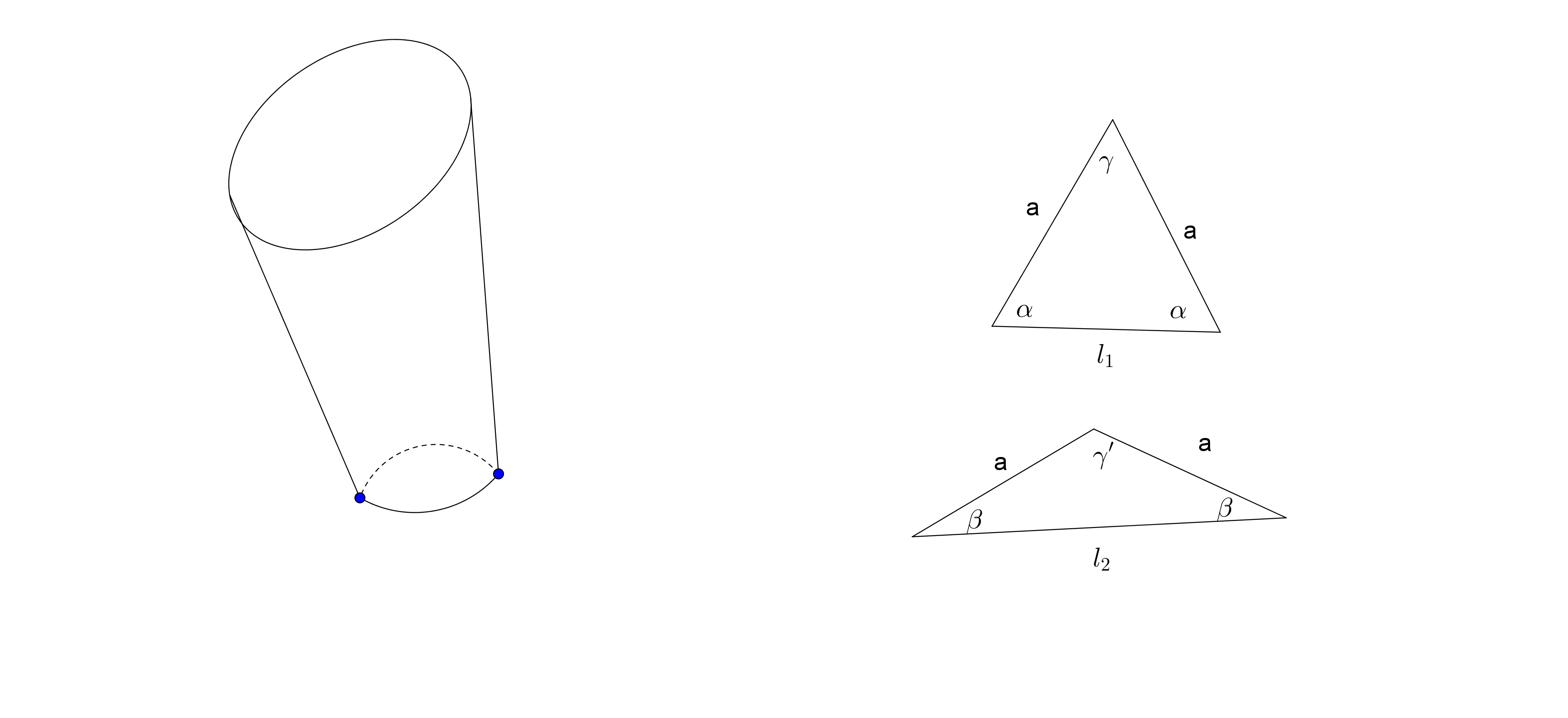}
	\end{center}
	\caption{Glue the isosceles triangles through the edges having length $a$ to get a flat disk with one singular interior point and two singular boundary points. Then glue this disk with $D^L$ to get a cone. }
	\label{cone2}
\end{figure}

\begin{lemma}
	\label{main-interior0}
	If $0 \leq \lvert \frak{K} \rvert < 2\pi$, then the puncture
on	any complete flat metric on $D^L$ is regular.

\begin{proof}
	We know that if $\frak{K}=0$ then $D^L$ is isometric to a half-cylinder. Thus the puncture is regular. 
	First assume that $0<\lvert \frak{K} \rvert < \pi$. By Theorem \ref{best-reduction}, this flat metric can be modified so that there is only one singular point on its boundary. Let $l$ be the length of its boundary. Take a isosceles triangle with angles $-\frak{K}, \frac{\pi + \kappa}{2}, \frac{\pi+\kappa}{2}$ with the length of the edge opposite to the vertex with angle $\frak{-K}$ is equal to $l$. If you glue the equal edges of this triangle and glue the resulting flat disk with $D^L$, you get the cone $C_{-\kappa}$. Hence the puncture is regular. See Figure \ref{cone1}.
	
	Assume that $\pi \leq \lvert \frak{K} \rvert < 2\pi$. Then we can modify $D^L$ so that it has two singular points of curvature $\frac{\frak{K}}{2}$. Call the singular points $b_1$ and $b_2$. Observe that the boundary has two components and these components connect $b_1$ and $b_2$. Let  $l_1$ and $l_2$ be length of these components. 
	Then we can find two isosceles triangles with the following properties:
	
	\begin{itemize}
		\item 
		First triangle has edges of length $l_1, a, a$ and the angles at its vertices are $\alpha,\alpha, \gamma$. Also the vertex with angle $\gamma$ is opposite the edge having length $l_1$.
		\item
		
		Second triangle has edges of length $l_2, a, a$ and the angles at its vertices are $\beta,\beta, \gamma$. Also the vertex with angle $\gamma'$ is opposite the edge having length $l_1$.
		\item
		$\gamma+\gamma'=-\frak{K}$.
	\end{itemize}
	Now glue the triangles along the edges having length $a$ to get a flat disk one singular interior point and two singular boundary points. Note that the angle at the singular interior point is $-\kappa$. If you appropriately glue this disk with $D^L$ along their boundaries, you get 
	the cone $C_{-\kappa}$. See Figure \ref{cone2}. Hence the puncture is regular.
	\end{proof} 

\end{lemma}

\subsection{Principal modifications }

Regarding Theorem \ref{best-reduction} , for each $\frak{K}<0$ so that  	$(n-1)\pi \leq \lvert \frak{K} \rvert < n \pi$, we will study modification-equivalence on the set

\begin{align*}
	\frak{C}(\frak{K},n)
	=\{\textup{Flat metrics on} \ D^L \ \textup{with n singular boundary points of curvature }\frac{\frak{K}}{n}\}/\textup{isometry}
\end{align*}

Note that two elements $\mu$ and $\eta$ are equivalent, $\mu \sim \eta$, if there is an orientation preserving isometry between them which respects the labeling of the vertices. 
Each element in $\frak{C}(\frak{K},n)$ is uniquely determined by the lengths of the boundary segments of the punctured disc. See \cite{saglam}[Theorem 1]. Let \\ $Ld(\mu)=(l[b_1,b_2],\dots ,l[b_n,b_{n+1}])$, where $l[b_i,b_{i+1}]$ is the length of the boundary segment joining $b_i$ and $b_{i+1}$. We call $Ld(\mu)$ as length data of $\mu$. Therefore  below map is a bijection:

$$ Ld: \frak{C}(\frak{K},n) \rightarrow \R_+^n$$
$$\mu \rightarrow Ld(\mu), $$
where $\R_+$ is the set of positive real numbers.

We will denote the set of equivalence classes of cone metrics on the disc by

$$\frak{M}(\frak{K})=\frak{M}(\frak{K},n):=\frak{C}(\frak{K},n)/\text{Modification}.$$ 
Now we define \textit{principal operations} which can be thought  as maps 
$\frak{C}(\frak{K},n) \rightarrow \frak{C}(\frak{K},n)$.

\begin{figure}
	\begin{center}
		\includegraphics[scale=0.45]{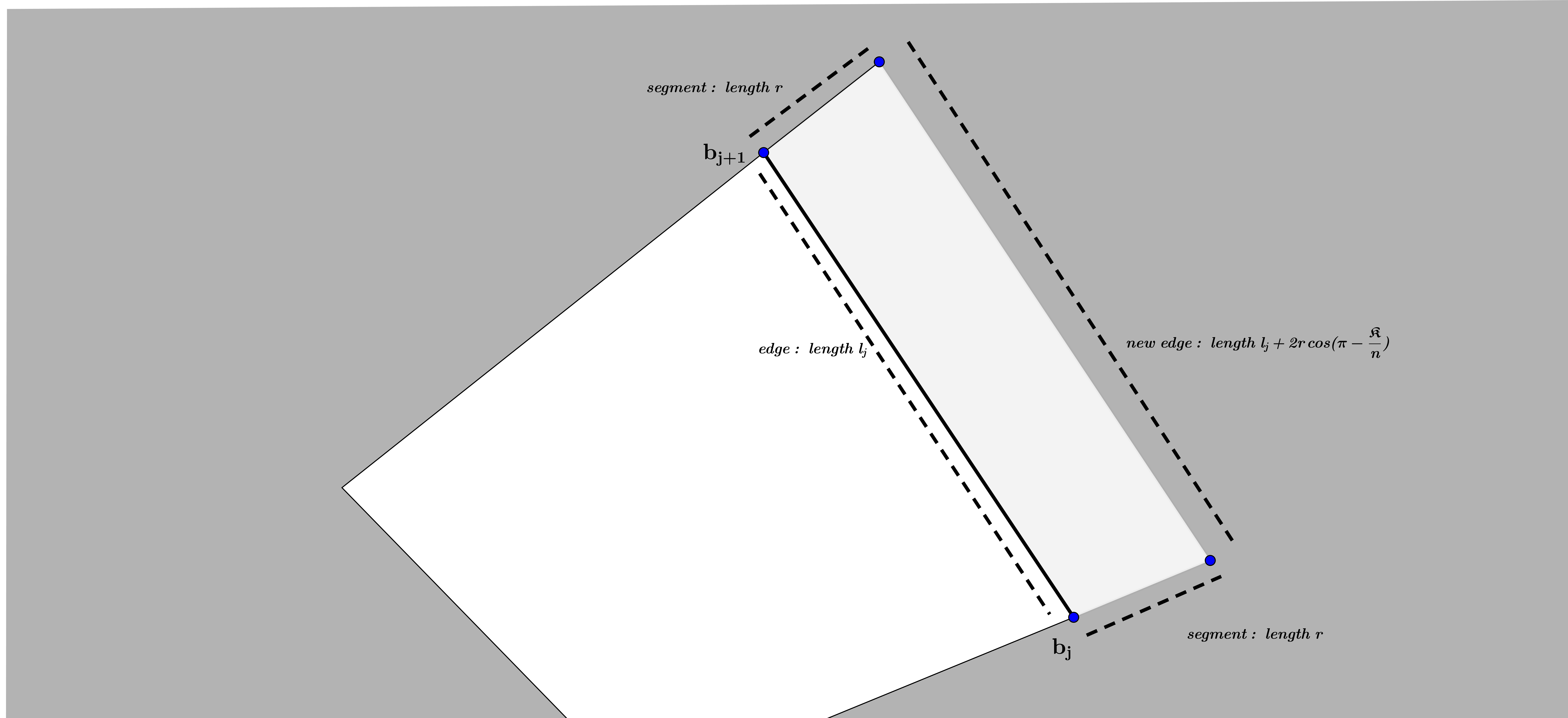}
	\end{center}
	\caption{An example of principal modifications : new cone metric is obtained by removing the quadrilateral in light grey. Note that $l_j=l[b_j,b_{j+1}]$   }
	\label{fig-principal-1}
\end{figure}

Let $n\geq 3$. Take an element in $\mu \in \frak{C}(\frak{K},n)$. Fix an  index $j \in \{1,2 \dots n \}$ and 
a non-negative real number $r$.  From the punctured disc (together with the metric $\mu$), subtract a quadrangle having angles

$$(\lvert \frac{\frak{K}}{n} \rvert, \pi - \lvert \frac{\frak{K}}{n} \rvert , \pi - \lvert \frac{\frak{K}}{n} \rvert,\lvert \frac{\frak{K}}{n} \rvert ),$$
and edge lengths 

$$(r,\ l[b_j,b_{j+1}]+ 2r \cos{(\pi - \frac{\frak{K}}{n})},\ r, \ l[b_j,b_{j+1}] ),$$

\noindent which has the segment $[b_j,b_{j+1}]$ as an edge. Since $\frac{\lvert \frak{K} \rvert}{n}+\frac{\lvert \frak{K} \rvert}{n}=2\frac{\lvert \frak{K} \rvert}{n}>\pi$, for each $r>0$ such a quadrangle exists. See Figure \ref{fig-principal-1}.
By this way, we obtain another element in $\frak{C}(\frak{K},n)$, denote this map by

$$\Theta_{j,r}: \frak{C}(\frak{K},n) \rightarrow \frak{C}(\frak{K},n).$$

\noindent From the Figure \ref{fig-principal-1}, description of this map in terms of length data is clear:

$$
\Theta_{j,r}: \R_+^n \rightarrow \R_+^n
$$
\begin{align}
	\Theta_{j,r}([l_1,\dots l_n])=[l_1,\dots,l_{j-2},l_{j-1}+r ,l_j+ 2r \cos(\pi - \frac{\frak{K}}{n}),l_{j+1}+r, l_{j+2}, \dots] \\
	=[l_1,\dots,l_n]+[0,\dots,0,\stackrel{j-1}{ \stackrel{\downarrow}{r}} ,2r \cos(\pi - \frac{\frak{K}}{n}),r, 0, \dots,0]
\end{align}

 Observe that above formulas imply the following:

\begin{enumerate}
	\item 
	$\Theta_{j,0}$ is identity map on $\frak{C}(\frak{K},n)$ (or on $\R^n$),
	\item
	$\Theta_{j,r}\circ \Theta_{j,r'}= \Theta_{j,r+r'}$,
	\item
	$\Theta_{j,r}\circ \Theta_{j',r'} = \Theta_{j',r'}\circ \Theta_{j,r}$.
\end{enumerate}
\begin{definition}
	We call semi-group generated by $\Theta_{j,r}$'s, either as maps on $\frak{C}(\frak{K},n)$ or $\R_+^n$,  the principal semigroup , and denote it as $\mathbb{T}=\mathbb{T}(\frak{K})$. 
\end{definition}



\begin{remark}
	\label{compact}
	Let $C$ be a compact set of $D^L$ and $\mu \in \frak{C}(\frak{K},n)$. There exists an element $\frak{T}$ in $\mathbb{T}(\frak{K})$ so that $C$ is a subset of the removed part of the once punctured sphere after the modification with respect to 
	$\frak{T}$.  
\end{remark}

\label{modification-equivalence}

\subsubsection{Circulant matrices}  We recall the basic properties of the circulant matrices.
\label{circulant-properties} A circulant matrix $\mathcal{C}$ is a $m\times m$ matrix obtained from 
one column vector $\mathbf{c}=[c_0,\dots,c_{m-1}]^T$  so that columns of $\mathcal{C}$ are determined 
by cyclic permutations of $\mathbf{c}$ as below:
\begin{align*}
	\mathcal{C}=
	\begin{bmatrix}
		c_0     & c_{m-1} & \dots  & c_{2} & c_{1}  \\
		c_{1} & c_0    & c_{m-1} &         & c_{2}  \\
		\vdots  & c_{1}& c_0    & \ddots  & \vdots   \\
		c_{m-2}  &        & \ddots & \ddots  & c_{m-1}   \\
		c_{m-1}  & c_{m-2} & \dots  & c_{1} & c_0 \\
	\end{bmatrix}.
\end{align*}

See \cite{circulant} for the  basic properties of circulant matrices. $f_{\mathcal{C}}(x)=c_0+c_1x+c_2x^2+\dots +c_{n-1}x^{n-1}$ is called associated polynomial of $\mathcal{C}$.

Let $\omega_j=\exp \left(\tfrac{2\pi \imath j}{n}\right)$ for each $j=0\dots,n-1$, where $\imath=\sqrt{-1}$. 

\begin{itemize}
	
	\item
	The set of eigenvalues of $\mathcal{C}$ is
	\begin{align}
		\{\lambda_j = c_0+c_{n-1} \omega_j + c_{n-2} \omega_j^2 + \ldots + c_{1} \omega_j^{n-1}: \ j=0, 1,\dots, n-1\}.
	\end{align}
	
	\item
	
	Determinant of $\mathcal{C}$, $\mathrm{det}(\mathcal{C})$ , is
	\begin{align}
		\label{determinant}
		\prod_{j=0}^{n-1} (c_0 + c_1 \omega_j + c_2 \omega_j^2 + \dots + c_{n-1}\omega_j^{n-1}).
	\end{align}
	
	\item
	
	Rank of $\mathcal{C}$ is the $n-d$ where $d$ is the degree of greatest common divisor of the polynomials $f_{\mathcal{C}}(x)$ and $x^{n}-1$.
	\item
	Eigenvector  with eigenvalue $\lambda_j$ is 
	$$v_j= [1,~ \omega_j,~ \omega_j^2,~ \ldots,~ \omega_j^{n-1}]^T,\quad j=0, 1,\ldots, n-1.$$
	Observe that a circulant matrix is diagonalizable and the eigenvectors of such a matrix do not depend on the coefficients $c_0, c_1, \dots, c_{n-1}$. 
\end{itemize}

\subsubsection{From the principal modifications to the circulant matrices}

Each $\Theta_{j,r}$ can be thought as a translation map $\R^n \to \R^n$. The vector space of translation maps can be identified with $\R^n$ and its canonical basis can be identified with the basis of $\R^n$ that consists of the vectors

$$e_1=(1,0,\dots,0)$$
$$e_2=(0,1, \dots,0)$$
$$\dots$$
$$e_n=(0,0,\dots,1)$$

\noindent With respect to this basis, $\Theta_{j,1}$ has coordinates 

$$[0,\dots,0,\stackrel{j-1}{ \stackrel{\downarrow}{1}} ,2 \cos(\pi - \frac{\frak{K}}{n}),1, 0, \dots,0]^T$$

\noindent Therefore the matrix that the coordinates of $\Theta_{1,1},\Theta_{1,1},\dots, \Theta_{n,1}$ form is a circulant matrix with $c_0=1,c_1=2 \cos(\pi - \frac{\frak{K}}{n}),c_2=1$ and $c_j=0$ if $j\neq 0,1,2$. Call this matrix $\mathcal{C}$. 

\begin{lemma}
	\label{determinant-key}
If $\frak{K}\neq -2\pi$, then $1+\cos{(\pi-\frac{\frak{K}}{n}})\omega_j+\omega_j^2\neq 0$.
	\begin{proof}
	Let  $\omega_j=\exp \left(\tfrac{2\pi \imath j}{n}\right)$. Assume that $(n-1)\pi \leq \lvert \frak{K}\rvert < n\pi$ and $n>3$. Then 
	$$(\frac{n-1}{n}+1)\pi \leq  \pi-\frac{\frak{K}}{n}< 2\pi$$
	$$\cos{\frac{\pi}{n}}\leq \cos{(\pi-\frac{\frak{K}}{n}}) < 1.$$
	In particular, $\frac{1}{2}<\cos{(\pi-\frac{\frak{K}}{n}})$. Now assume that $1+2\cos{(\pi-\frac{\frak{K}}{n}})\omega_j+\omega_j^2=0$. It follows that $\lvert 1+\omega_j^2 \rvert= 2\cos{(\pi-\frac{\frak{K}}{n}})>2$, which is impossible. Hence $1+2\cos{(\pi-\frac{\frak{K}}{n}})\omega_j+\omega_j^2\neq 0$ for all $j$. A similar argument shows that $1+2\cos{(\pi-\frac{\frak{K}}{n}})\omega_j+\omega_j^2\neq 0$ when $n=3$ and $\kappa\neq -2\pi$.
	\end{proof}
\end{lemma}
\begin{proposition}
	\label{main-interior1}
	If $\kappa< - 2\pi$, then $\frak{M}(\kappa)$ consists of one single point, that
	is, there are no irregular punctures on $D^L$ when the total curvature on the boundary is not equal to $-2\pi$.
	\begin{proof}
We will show that $\frak{C}(\frak{K},n)/\mathbb{T}(\frak{K})$ consists of a single point. Clearly this implies that $\frak{M}(\kappa)$ consists of 	one point. Lemma \ref{determinant-key} implies that det$(\mathcal{C})\neq 0$. Therefore the group generated by $\Theta_{1,1},\dots,\Theta_{n,1}$ is the full group of translations of $\R^n$. So this group has one orbit. It follows that $$\frak{C}(\frak{K},n)/\mathbb{T}(\frak{K})=\R_+^n/\mathbb{T}(\frak{K})$$
has only one point.
	\end{proof}
\end{proposition}

\noindent Now we consider the case $\frak{K}=-2\pi$. 

\begin{proposition}
	\label{main-interior2}
There is a bijection between $\frak{M}(-2\pi)$ and $\R^2_+$.
\begin{proof}
	Consider the map  $\R^2_+\to \frak{M}(-2\pi)$ sending $(\alpha,\beta)$
	to the modification-equivalence class of the metric in $\frak{C}(-2\pi,3)$ having boundary segments of length $1,\alpha, \beta$. Let us denote the modification-equivalence class of the metric on $D^L$ with boundary segment of length $a,b,c$ by $[a,b,c]$. Since $\cos (\pi-\frac{2\pi}{3})=\frac{1}{2}$, it follows that $\Theta_{1,r}$ sends $[a,b,c]$ to $[a+r,b+r,c+r]$. Therefore it is easy to see that this map is surjective. Now assume that $[1,b,c]$ and $[1,b',c']$ are modification equivalent. It follows that there is an orientation preserving isometry sending the flat disk $D^L$ whose boundary segments have length $1+r,b+r,c+r$ to the flat disk $D^L$ whose boundary segments have length $1+r',b'+r',c'+r'$, where $r,r' \geq 0$. It follows that $r=r'$, $b=b'$ and $c=c'$. So the map is injective.
\end{proof}
\end{proposition}

Now we collect the results in Lemma \ref{main-interior0}, Proposition \ref{main-interior1} and Proposition \ref{main-interior2} in a single theorem.

\begin{theorem}
	\label{main-interior}
	\begin{enumerate}
		\item
	If $\frak{K}\neq-2\pi$, then $\frak{M(\frak{K})}$ consists of one singular point. In other words there are no complete flat metrics on $D^L$ so that the puncture is irregular when curvature at the puncture is not equal to $4\pi$.
\item 
$\frak{M}(-2\pi)\equiv \R^2_+ $. 

\end{enumerate}
	 
\end{theorem}

\section*{Acknowledgements}
I am really grateful to Ahmet Refah Torun  for his remarks and corrections. This work is supported by Research Fund of Adana Alparslan T\"{u}rke\c{s} Science and Technology University. Project Number: 18119001.

\end{document}